\documentstyle{article}
\input amssym.def
\newcommand{\unit}{{\bf \cal I}}
\newcommand{\qs}{\ast}
\begin{document}

\begin{center} {\Large\bf The Geometry of Deformed Boson Algebras\\}
\vspace{5mm}
{\large \bf P. Crehan\footnote{This work is dedicated to John Kennedy
my former teacher who died of cancer 20th of January 1995}} \\
{\it Dept. of Mathematics, Kyoto University\\
pat@kusm.kyoto-u.ac.jp \\}
{\large \bf T. G. Ho}\\
{\it Research Institute for Mathematical Sciences, Kyoto\\
 tingguo@kurims.kyoto-u.ac.jp }
\end{center}

\begin{abstract}
Phase-space realisations of
an infinite parameter family of quantum deformations of the boson algebra
in which the $q$-- and the $qp$--deformed
algebras arise as special cases are studied.
Quantum and classical models for the corresponding deformed
oscillators are provided. The deformation parameters are identified with
coefficients of non-linear terms in the normal forms expansion
of a family of classical Hamiltonian systems.
These quantum deformations are trivial in the sense
that they correspond to non-unitary transformations of the Weyl algebra.
They are non-trivial in the
sense that the deformed commutators consistently
quantise a class of non-canonical classical Poisson structures.
\end{abstract}
\newpage
\section{Introduction}
It has been pointed out that the quantum group
is not an intrinsically quantum concept. Flato and Lu
have provided a purely classical realisation of $U({\goth sl}_{q}(2))$ in the
canonical Poisson algebra of functions on ${\Bbb R}^2$ \cite{flato_lu}, while
Chang, Chen and Guo have done the same on ${\Bbb R}^4$
\cite{chang_chen_guo}. In particular they have
constructed classical Poisson brackets which depend on
the $q$-deformation parameter, such that the fundamental Poisson bracket
relations are those of $U({\goth sl}_{q}(2))$. They conclude
that although $q$-deformation is related to $*$-deformation,
it is not the same as quantisation. Their results demonstrate
a need to clarify the difference between quantum-deformation and
quantisation, and to understand the essential role
of the deformation parameter. In a recent study of the
$q$-boson, Man'ko, Marmo, Solimeno and
Zaccaria propose that the role of the deformation parameter is to
introduce non-linearity in the deformed oscillator \cite{manko}.

{}From a different point of view the Weyl algebra is known to be rigid
\cite{weyl:rigidity}. It
does not possess a non-trivial associative deformation in the sense that
any associative deformation of the Weyl algebra
is formally isomorphic to the Weyl algebra itself.
This applies in particular to quantum deformations of
the Weyl algebra related to the $q$-- and the $qp$--bosons.
The isomorphisms which relate the deformed algebras
to the Weyl algebra cannot be unitary since
they do not preserve the fundamental commutation relations.
The real problem therefore is to understand these non-unitary
equivalence transformations. We would like to know if quantum
deformation has a natural geometrical meaning.

In {\S}~\ref{boson_algebras} we introduce families of
deformations of the Weyl algebra and of the Heisenberg Lie algebra
which depend on an infinite number of independent parameters.
In {\S}~\ref{phase-space_quantum_mechanics} we introduce the tools we will
need from phase-space quantum mechanics. We extend these to the
deformed phase-space in {\S}~\ref{extension}.
The corresponding Moyal algebras are studied using
bimodular bases of deformed oscillator eigenstates and
bases of anti-standard ordered monomials in {\S}~\ref{bimodular_basis} and
{\S}~\ref{aso_basis} respectively. In the eigenstate bases it is easy to see
the isomorphism between the deformed algebras. In {\S}~\ref{infinite_matrices}
we demonstrate an isomorphism between the deformed algebras
and the algebra of infinite matrices. Quantum models of the
deformed algebras and their associated oscillators are
dealt with in {\S}~\ref{quantum_models}.
In {\S}~\ref{nonlin} we show how phase-space representations of
the deformed algebras are related by
non-linear non-unitary transformations. We discuss the provision of
co-product and Hopf-algebra structures in {\S}~\ref{coproduct}.

In {\S}~\ref{classical_limit} we consider the classical limit of the
deformed quantum dynamics. The deformed commutators contract
to non-canonical Poisson brackets. We identify the deformed Weyl algebras
with non-linear deformations of the harmonic oscillator, and
in {\S}~\ref{normal_forms} we reveal a relationship with the classical
method of normal forms.

The overall picture which emerges is that quantum
deformations of the Weyl algebra
are trivial in the sense that they correspond to non-linear
non-unitary equivalence transformations of $U(W(1))$.
The parameters which arise in quantum deformation of the Weyl algebra
have a geometrical meaning in that they parameterise a family of
non-canonical transformations of classical phase-space.
The quantum deformations are non-trivial in the sense that they
quantise a family of non-canonical classical Poisson brackets.
In {\S}~\ref{quantisation}
we reinterpret the algebras introduced in \S~\ref{boson_algebras}
as quantum algebras which consistently quantise
the infinite number of Hamiltonian
structures of the simple harmonic oscillator.
This indicates a number  of directions in which certain
applications of quantum groups can naturally be extended.

\section{Deformed boson algebras}\label{boson_algebras}
The defining relations for the $q$-deformed boson algebra have been
given in many different forms \cite{kulish_2}, \cite{mcfarlane}.
Taking $A$, and $A^+$ to be the $q$-deformed
annihilation and creation operators, the number operator $N$ always
satisfies

\begin{eqnarray}
{[} N, A ] =  -A &&{[} N, A^+ ] = A^+        \label{N_brackets}
\end{eqnarray}

\noindent
while $A$ and $A^+$ are related as follows

\begin{eqnarray}
A A^+-q A^+ A&=&q^{-N}. \label{q_boson}
\end{eqnarray}

\noindent
A generalisation known as the $qp$-boson
\cite{chakrabarti_jagannathan}, \cite{kibler} satisfies

\begin{eqnarray}
A  A^+- q A^+  A&=&{q^{N+1}-p^{N+1} \over q-p}- {q^N-p^N \over q-p}.
                                                         \label{qp_boson}
\end{eqnarray}

\noindent
A further generalisation satisfies

\begin{eqnarray}
A  A^+-q A^+  A&=&P_{q}(N) \label{won_boson}
\end{eqnarray}

\noindent
where $P_{q}(N)$ is a polynomial of finite order in the number
operator $N$ whose coefficients functionally depend on a single
deformation parameter $q$ \cite{won}.

The relations (\ref{N_brackets}) along with
(\ref{q_boson}), (\ref{qp_boson}) or (\ref{won_boson})
constitute re-ordering rules or equivalence relations for the
infinite polynomial ring generated by $A$, $A^+$, $N$ and $1$.
They are deformations  of  the standard
boson algebra in the sense that in the limit $q,p\to 1$,
(\ref{q_boson}), (\ref{qp_boson}) and (\ref{won_boson}) reduce to

\begin{eqnarray}
A  A^+ - A^+ A & = & 1. \label{boson}
\end{eqnarray}

In any given
representation these re-ordering rules can be replaced by
equivalent commutation relations.
In the case of the $q$-deformed boson algebra for example,
(\ref{q_boson}) can be replaced by
one of the following commutation relations, one for each of the three
families of representations of $A A^+ - q A^+ A = q^{-N}$
provided by Rideau \cite{rideau}.

\begin{eqnarray}
{[}A,A^+]&=& q^{-\nu_0} ([N+1]_q - [N]_q)                       \nonumber \\
{[}A,A^+]&=& q^{-\nu_0} ([N+1]_q - [N]_q) + \lambda_0 (q-1) q^N \nonumber \\
{[}A,A^+]&=& - q^{-\nu_0}(1+q)^{-1} q^{-N}.        \label{commutator_form}
\end{eqnarray}

\noindent
We have modified our notion slightly so that $N$ satisfies
$N|n\rangle = n |n\rangle$ on the relevant Fock-space.
The advantage to be gained by rewriting (\ref{q_boson}) in this
way is that the equivalence relations now have a natural dynamical
interpretation. On the basis of (\ref{commutator_form})
we propose that any representation of a deformed Weyl-algebra
has an equivalent formulation in terms of commutators. {}From now on we
will consider algebras which consist of the infinite ring
of polynomials generated by $A$, $A^+$, $N$ and $1$ modulo the
following relations

\begin{eqnarray}
{[} A, N ]     &=& A       \nonumber                           \\
{[} N, A^+ ]   &=& A^+     \nonumber                           \\
{[} A, A^+ ] &=& f_{\bar{q}}(N).     \label{wqbar:boson}
\end{eqnarray}

\noindent
We could have replaced the commutator of $A$ and $A^+$ with
something of the form
$A A^+ - q A^+ A = f_{\bar{q}}(N)$. Doing so makes no
difference to our results and the form that we have
chosen has the advantage of being easier to interpret dynamically.
For simplicty we assume that $f_{\bar{q}}$ is an entire function
of the complex plane such that $f_{\bar{q}}({\Bbb R}^+)\ge 0$.
Everything we do will be valid
for functions which are real analytic
and positive on $[0,\infty )$. In fact $f_{\bar{q}}(x)$
is not so  much a function as an equivalence class of functions
defined by the sequence of values taken on the positive integers.
In a neighbourhood of the origin we identify $f_{\bar{q}}$
with a power series expansion
$ f_{\bar{q}}(N) = q_0 + q_1 N + q_2 N^2 +\ldots$ whose coefficients
satisfy $q_i\in {\Bbb R}$ for all $0\le i\in {\Bbb Z}$ and
constitute independent deformation
parameters. We denote this category of deformed algebras by $\Bbb W$.

The Heisenberg algebra ${\goth h}(1)$ is closely related to the
Weyl-algebra. It is
the Lie algebra with three generators
$A$, $A^+$ and $E$ whose only non-vanishing bracket is given by
$[A,A^+ ] = \hbar E$. {}From a physical point of view the Weyl-algebra
corresponds to a phase-space realisation of ${\goth h}(1)$, it is
the quotient of ${\goth h}(1)$ with the two-sided ideal $E=1$, it is not a Lie
algebra, and according to one authority is not known
to possess a Hopf-algebra structure \cite{kulish_1}.
The Heisenberg Lie algebra does have
a Hopf-algebra structure, and
some authors have shown that it is instructive to consider In\"on\"ou-Wigner
contractions {}from $q$-deformed semi-simple quantum groups such as
$U({\goth {su}}_q(2))$ to the $q$-deformed Heisenberg
group $U({\goth h}_q(1))$ \cite{celeghini}.
There are advantages to be gained {}from considering the
deformation of boson algebras
{}from the point of view of the Heisenberg Lie algebra, but there
is no unique way in which
to extend Weyl-type algebras satisfying
relations such as (\ref{q_boson}) or (\ref{commutator_form}) to
Heisenberg-type algebras. For example one could replace (\ref{q_boson}) with
$A A^+-q^E A^+ A=E q^{-N}$, with  $A A^+ - q E A^+ A = E q^{-N}$
or with any one of an infinite number of other choices.
In any case we will keep in mind that it is possible to
consider a family of deformed algebras of Heisenberg type
generated by $A$, $A^+$, $E$ and $N$, whose non-vanishing
commutation relations can be given in the following form

\begin{eqnarray}
{[} A, N ]     &=& E A       \nonumber                           \\
{[} N, A^+ ]   &=& E A^+     \nonumber                           \\
{[} A, A^+ ] &=& f_{\bar{q}}(N,E).     \label{hqbar:boson}
\end{eqnarray}

\noindent
Various restrictions should apply to $f_{\bar{q}}(x,y)$ so that
it it is analytic and positive on ${\Bbb R}^+\times {\Bbb R}^+$.
As before the coefficients of its power series expansion about
the origin can be thought of as independent deformation parameters.
The infinite polynomial ring generated by  $A$, $A^+$, $E$ and $N$
modulo (\ref{hqbar:boson}) is a deformation of $U({\goth h}(1))$ in the sense
that for suitable values of these coefficients
we recover the commutation relations of the Heisenberg
algebra with $N=A^+ A$. We denote these
deformations of the Heisenberg Lie algebra by ${\goth H}_{\bar{q}}$.

In the next section we summarise relevant results
{}from the phase-space theory of the boson oscillator,
which we adapt to study $W_{\bar{q}}\in \Bbb W$,
the related  ${\goth H}_{\bar{q}}$ and the deformed phase-space dynamics of
their associated classical and quantum oscillators.

\section{Deformed phase-space dynamics}\label{phase-space_quantum_mechanics}
In the phase-space quantum mechanical formalism the space of states has a
bimodular structure like that of the algebra of observables.
By considering the eigenvalue problem of the harmonic oscillator
it is possible to construct a basis for the space of states which
consists of left-right eigenstates of the harmonic oscillator Hamiltonian.
Eventually the algebra of observables can be thought of as lying in the
closure of finite linear combinations of states with respect to the
appropriate topology. We will use an extension of results such as these
in developing the theory of deformed bosonic algebras. We start by
describing in more detail the case of the harmonic oscillator.

In this formalism the vacuum state is required to satisfy

\begin{eqnarray}
a*\Omega_{00} = &0& = \Omega_{00}*a^+,        \label{ho:vacuum_post}
\end{eqnarray}

\noindent
where $f*g$ stands for the Moyal product of $f$ and $g$.
The solution to this  equation is unique up to a constant
and when properly normalised satisfies
$\Omega_{00}*\Omega_{00}=\Omega_{00}$. The left-right eigenstates
of the harmonic oscillator are given by

\begin{eqnarray}
\Omega_{nm}&=&(n!\ m! \ \hbar^{n+m})^{-1/2} {a^+}^n * \Omega_{00}*a^m
                					\label{ho:states}
\end{eqnarray}

\noindent
and have the following remarkable properties

\begin{eqnarray}
\overline{\Omega_{n m}(q,p)}&=&\Omega_{m n}(q,p)                \nonumber \\
\Omega_{n m}(q,p)*\Omega_{n^\prime\ m^\prime}(q,p)&=&
\delta_{m n^\prime}\ \Omega_{n m^\prime}(q,p) \label{ho:product}
\end{eqnarray}

\noindent
as well as
\begin{eqnarray}
(2 \pi \hbar)^{-1} \int_{{\Bbb R}^2} \overline{\Omega_{n m}(q,p)}\
\Omega_{n^\prime
m^{\prime}}(q,p)\ \ dq\wedge dp &=& \delta_{n n^{\prime}}\
\delta_{m m^{\prime}}                     \nonumber       \\
(2 \pi \hbar)^{-1} \int_{{\Bbb R}^2} \Omega_{n m}(q,p)\ dq\wedge dp &=
   & \delta_{n m}.                          \label{ip:boson}
\end{eqnarray}

\noindent
These functions were introduced by Hansen
who showed that they provide a basis for the Schwartz space $S({\Bbb R}^2)$,
the Hilbert  space of
quantum states on phase-space $L^2({\Bbb R}^2, dq\wedge dp)$
and for $S^{\prime}({\Bbb R}^2)$ the space of tempered distributions of
$S({\Bbb R}^2)$
\cite{hansen}. It is customary to study $W$ as a quotient of
the infinite polynomial
ring generated by $a$ and $a^+$. Instead of using a basis of
Weyl-ordered monomials, we found it more convenient to work with
Anti-Standard Ordered Monomials ASOM's. In this basis
the intertwiner of $W$ is given by

\begin{eqnarray}
f(a,a^+) \bullet g(a,a^+) &=& \exp\left(\hbar\ \partial^{(1)}_{a^+}\cdot
                  \partial^{(2)}_{a}\right) \
                          f^{(1)}(a,a^+)\ g^{(2)}(a,a^+)
                                           \label{aso:prod}
\end{eqnarray}

\noindent
where $f$ and $g$ are the ASO representations of elements of $W$.
An identity distribution
in $S^{\prime}({\Bbb R}^2)$ which satisfies

\begin{eqnarray}
\unit*f(q,p)&=&f(q,p) \nonumber \\
\int_{{\Bbb R}^2} \unit*g(q,p)\ dq\wedge dp&=& \int_{{\Bbb R}^2} g(q,p)\
dq\wedge dp
			\label{unit:def}
\end{eqnarray}

\noindent
for all $f(q,p)\in S^{\prime}({\Bbb R}^2)$, and $g(q,p)\in S({\Bbb R}^2)$, and
which is given by

\begin{eqnarray}
\unit &=& \sum_{n=0}^{\infty} \Omega_{n n}(q,p) \label{unit:rep}
\end{eqnarray}

\noindent
was used to construct the linear transformation {}from the
ASO basis for $W$ to the $\Omega_{n m}$ basis.
The following expression for the vacuum state

$$\Omega_{00}= 2 \exp(-{2\over \hbar} a a^+ )=
          \sum_{n=0}^{\infty} { (-1)^n\over \hbar^n\ n!} {a^+}^n * a^n$$

\noindent
was used to construct the
inverse transformation {}from the left-right eigenstate basis
back to the basis of ASOMs.

By providing an algorithm for moving between them,
these results establish the equivalence of the ASO monomial basis
and the left-right eigenvector basis for the Moyal algebra of observables.
We now turn to the infinite parameter family of deformations
of the Weyl algebra (\ref{wqbar:boson}), which we study in the monomial
and oscillator eigenstate bases.

\subsection{Extension to the quantum deformation $W_{\bar{q}}$}
\label{extension}
We assume the existence of a deformed
$\qs$-product which satisfies (\ref{wqbar:boson}), and
postulate the existence of a vacuum state $\Omega^{\bar{q}}_{00}$
with the following reasonable properties

\begin{eqnarray}
\begin{array}{rcl}
A\qs \Omega^{\bar{q}}_{00}&=&0 \\
N\qs \Omega^{\bar{q}}_{00}&=&0 \\
\end{array}
&&
\begin{array}{rcl}
\Omega^{\bar{q}}_{00} \qs A^+&=&0 \\
\Omega^{\bar{q}}_{00} \qs N  &=&0 \\
\end{array}
   \nonumber  \\
&& \nonumber  \\
\Omega^{\bar{q}}_{00}\qs \Omega^{\bar{q}}_{00}&=&\Omega^{\bar{q}}_{00}.
                                                 \label{qbar:vacuum}
\end{eqnarray}

\noindent
When considering ${\goth H}_{\bar{q}}$ it will also be necessary to add
$\Omega^{\bar{q}}_{00} \qs E = E\qs \Omega^{\bar{q}}_{00} =
\Omega^{\bar{q}}_{00}$. {}From (\ref{wqbar:boson}) we define the following
real functions on the strictly positive integers

\begin{eqnarray*}
F_{\bar{q}}(n)&=& \sum_{i=1}^{n} f_{\bar{q}}(i-1)  \\
F_{\bar{q}}!(n)&=& \prod_{i=1}^{n} F_{\bar{q}}(i),
\end{eqnarray*}

\noindent
and extend them  to $n=0$ by putting $F_{\bar{q}}(0)=0$ and
$F_{\bar{q}}!(0)=1$. Due to the restrictions imposed on $f_{\bar{q}}$
in \S~\ref{boson_algebras}, $F_{\bar{q}}(n)$ is positive and monotonic
increasing function on $0\le n \in {\Bbb Z}$.
 {}From the re-ordering relations (\ref{wqbar:boson}),
and the vacuum postulate (\ref{qbar:vacuum}), one can derive the following
identities

\begin{eqnarray}
A\qs {A^+}^n \qs \Omega^{\bar{q}}_{00}&=&
     F_{\bar{q}}(n)\ {A^+}^{n-1}\qs \Omega^{\bar{q}}_{00},     \nonumber\\
\Omega^{\bar{q}}_{00}\qs  {A^n} \qs {A^+}^m\qs \Omega^{\bar{q}}_{00}&=&
                                       \delta_{n m}\ F_{\bar{q}}!(n)\
         \Omega^{\bar{q}}_{00}.              \label{identities}
\end{eqnarray}

The $\bar{q}$-states are defined as follows

\begin{eqnarray}
\Omega^{\bar{q}}_{nm}&=&
                 \left( F_{\bar{q}}!(n) F_{\bar{q}}!(m) \right)^{-1/2}\
        {A^+}^n\qs \Omega^{\bar{q}}_{00}\qs A^m.\label{qbar:states}
\end{eqnarray}

\noindent
{}From (\ref{identities}) it is easily seen that these are
left-right eigenvectors of the $\bar{q}$-deformed oscillator Hamiltonian

\begin{eqnarray}
H &=&{1\over 2}\left( A^+\qs A + A \qs A^+ \right)  \label{qbar:hamiltonian}
\end{eqnarray}

\noindent
since

\begin{eqnarray}
H\qs \Omega^{\bar{q}}_{n m}&=&E^{\bar{q}}_n\ \Omega^{\bar{q}}_{n m}\nonumber\\
\Omega^{\bar{q}}_{nm}\qs H&=&E^{\bar{q}}_m\ \Omega^{\bar{q}}_{n m}\nonumber \\
E^{\bar{q}}_n&=& {1\over 2}\left( F_{\bar{q}}(n+1) + F_{\bar{q}}(n)\right).
                                                             \label{q_spectrum}
\end{eqnarray}

\noindent
The standard boson corresponds to $f_{\bar{q}}(n)=\hbar$, in which case
$F_{\bar{q}}(n)=n\hbar$, $F_{\bar{q}}!(n)=n!\hbar^n$, and the
eigenvalues $E^{\bar{q}}_n$ are the familiar $\hbar(n+1/2)$.
Later we will use the following relations which can be derived {}from the
identities (\ref{identities})
and the definitions (\ref{qbar:states})

\begin{eqnarray}
\begin{array}{rcl}
E\qs \Omega^{\bar{q}}_{n m}  &=&\Omega^{\bar{q}}_{n  m}\\
A\qs \Omega^{\bar{q}}_{n m}  &=&F_{\bar{q}} (n)^{1/2}\
                                                 \Omega^{\bar{q}}_{n-1\ m}\\
A^+\qs \Omega^{\bar{q}}_{n m}&=&F_{\bar{q}} (n+1)^{1/2}\
                                         \Omega^{\bar{q}}_{n+1\ m}\\
N\qs \Omega^{\bar{q}}_{n m}  &=&n\ \Omega^{\bar{q}}_{n m}\\
\end{array}
&\ &
\nonumber\\
\begin{array}{rcl}
\Omega^{\bar{q}}_{n m}\qs E &=&\Omega^{\bar{q}}_{n m}\\
\Omega^{\bar{q}}_{n m}\qs A&=&F_{\bar{q}} (m+1)^{1/2}\
                                      \Omega^{\bar{q}}_{n m+1}\\
\Omega^{\bar{q}}_{n m}\qs A^+&=&F_{\bar{q}}(m)^{1/2}\
                                        \Omega^{\bar{q}}_{n m-1}\\
\Omega^{\bar{q}}_{n m}\qs N &=& m\ \Omega^{\bar{q}}_{n m}.\\
\end{array}
                                                    \label{q_ladder}
\end{eqnarray}

Using the projection property of the vacuum (\ref{qbar:vacuum})
and the identities (\ref{identities}), one can now show that in this
basis the deformed Moyal products are given by

\begin{eqnarray}
\overline{\Omega^{\bar{q}}_{n m}}&=& \Omega^{\bar{q}}_{m n} \nonumber     \\
\Omega^{\bar{q}}_{n m}\qs \Omega^{\bar{q}}_{n^\prime m^\prime}&=&
                      \delta_{m n^\prime}\
\Omega^{\bar{q}}_{n m^\prime}.               \label{algebra_of_states}
\end{eqnarray}

\noindent
Just as the Weyl algebra   (\ref{boson}) was  replaced   by
(\ref{ho:product}), the defining relations of the   deformed
Weyl  algebras (\ref{wqbar:boson})  have been replaced  by
(\ref{algebra_of_states}).  In these oscillator eigenstate bases
the deformed Moyal products do not depend on values of the
deformation parameters. The expected isomorphism between the  Weyl
algebra and its deformations becomes
immediately  apparent. In the case of phase-space     realisations   of
(\ref{wqbar:boson})    which         satisfy  the vacuum     postulate
(\ref{qbar:vacuum}), this  turns out to be simply an isomorphism with  the
algebra of infinite matrices.

\noindent
By analogy with (\ref{ip:boson}) the natural inner product is given by

\begin{eqnarray}
\langle \Omega^{\bar{q}}_{n m}, \Omega^{\bar{q}}_{n^\prime m^\prime}
\rangle &=& \delta_{n n^\prime}\ \delta_{m m^\prime}. \label{ip:qboson}
\end{eqnarray}

\noindent
{}From Hansen's work \cite{hansen}, we know that the closure of the
space of finite linear combinations of $\Omega^{\bar{q}}_{nm}$ with respect to
appropriate topologies will provide the deformed versions
of $S({\Bbb R}^2)$, $L^2({\Bbb R}^2, dq\wedge dp)$ and
$S^\prime({\Bbb R}^2)$. These
can now be identified with the undeformed spaces.

In view  of  the rigidity  of  the Weyl  algebra  the   equivalence of
phase-space  realisations  of    the  infinite  parameter  family   of
$\bar{q}$-boson algebras  to  that of the standard boson algebra  does not
come  as a surprise.  However this  equivalence  does not preserve the
defining relations (\ref{wqbar:boson}), and so is different {}from
the unitary   equivalence  of quantum   mechanics. In {\S}~\ref{nonlin}
we apply elements {}from a theory  of non-linear representations of
Lie  algebras to   understand  this equivalence
better.  In the  next section we establish  the existence of  concrete
representations of   $W_{\bar{q}}$,  demonstrate the equivalence of
oscillator
eigenstate bases to  bases of ASO  monomials and  we provide classical
and quantum models for the deformed oscillators.

\section{Realisations and Representations of ${\goth H}_{\bar{q}}$ and
$W_{\bar{q}}$}
\subsection{Realisation of ${\goth H}_{\bar{q}}$ and
$W_{\bar{q}}$ as $\{ \Omega_{\bar{q}}, * \}$.} \label{bimodular_basis}

We denote by $\Omega_{\bar{q}}$ the set of formal linear combinations
of  $\Omega^{\bar{q}}_{n m}$, and by $\Omega$ the set of formal
linear combinations of left-right eigenstates of the standard boson.
Equipped with a $*$--product obeying (\ref{algebra_of_states}) we know {}from
Hansen's work that suitably restricted subsets of $\{ \Omega_{\bar{q}}, *\}$
are Moyal algebras \cite{hansen}.
We will look at a family of projections {}from  ${\goth H}_{\bar{q}}$
($W_{\bar{q}}$)
to $\{ \Omega_{\bar{q}}, * \}$ which are ${\goth H}_{\bar{q}}$
($W_{\bar{q}}$) algebra homomorphisms.

In particular we consider $\Sigma_{\bar{q}}$
{}from ${\goth H}_{\bar{q}}$ to $\{\Omega_{\bar{q}}, *\}$ based on $\unit$
a resolution of the identity on  $\{\Omega_{\bar{q}}, *\}$. $\unit$
satisfies $\unit \qs \omega_{\bar{q}} = \omega_{\bar{q}} \qs \unit =
\omega_{\bar{q}}$ for all $\omega_{\bar{q}} \in \Omega_{\bar{q}}$,
and is given by

$$ \unit = \sum_{n=0}^\infty \Omega^{\bar{q}}_{n n}.$$

\noindent
The ${\goth H}_{\bar{q}}$--homomorphism $\Sigma_{\bar{q}}$ is defined by

\begin{eqnarray}
\Sigma_{\bar{q}} ( u ) &=& \unit \qs u \qs \unit \label{bigsigma:barq}
\end{eqnarray}

\noindent
for all $u \in {\goth H}_{\bar{q}}$. Its action on
$u \in {\goth H}_{\bar{q}}$ can be computed explicitly using the
ladder relations (\ref{q_ladder}). These relations
can be used to show that for $u_1$ and $u_2$ in ${\goth H}_{\bar{q}}$ we have

$$\Sigma_{\bar{q}}(u_1) \qs \Sigma_{\bar{q}}(u_2)
 = \Sigma_{\bar{q}}(u_1 \qs u_2),$$

\noindent so $\Sigma_{\bar{q}}$ is a deformed algebra homomorphism {}from
${\goth H}_{\bar{q}}$ into $\{\Omega_{\bar{q}},*\}$.
Given a realisation of the $\qs$-product (\ref{wqbar:boson})
which admits a vacuum state satisfying the requirements
(\ref{qbar:vacuum}), $\Sigma_{\bar{q}}$ provides a
representation of ${\goth H}_{\bar{q}}$, and therefore of $W_{\bar{q}}$
in $\Omega_{\bar{q}}$ acting on $\Omega_{\bar{q}}$. We have

\begin{eqnarray*}
E&\to&\sum_{n=0}^{\infty} \Omega^{\bar{q}}_{n n}(q,p)\\
A^+&\to&\sum_{n=0}^\infty F_{\bar{q}}(n+1)^{1/2}\
                 \Omega^{\bar{q}}_{n+1\ n}(q,p)\\
A&\to&\sum_{n=0}^\infty F_{\bar{q}}(n+1)^{1/2}\
                \Omega^{\bar{q}}_{n\ n+1}(q,p)\\
N&\to&\sum_{n=0}^\infty n\ \Omega^{\bar{q}}_{n n}(q,p)\\
H&\to& {1\over 2} \sum_{n=0}^\infty (F_{\bar{q}}(n+1)+F_{\bar{q}}(n))\
  \Omega^{\bar{q}}_{n n}(q,p).
\end{eqnarray*}

\noindent
These infinite sums are meaningful in the sense of distributions on
$S^{\prime}({\Bbb R}^2)$. The action of these operators on the Hilbert
space spanned by $\Omega_{\bar{q}}$ can be computed directly {}from
(\ref{algebra_of_states}).

\subsection{Isomorphism of $\{\Omega_{\bar{q}}, *\}$ and
$\{W^{AS}_{\bar{q}}, \bullet\}$.}\label{aso_basis}

We denote by $W^{AS}_{\bar{q}} \subseteq {\goth H}_{\bar{q}}$ the linear
space spanned by ${A^+}^n \qs A^m $ for $n,m \ge 0$, the
ASO monomial functions of co-ordinates on a $\bar{q}$-deformed
phase-space. The restriction of $\Sigma_{\bar{q}}$
to $W^{AS}_{\bar{q}}$ will be denoted $\sigma_{\bar{q}}$
and $\sigma$ refers to its action on $W^{AS}$ the space of formal linear
combinations of ASO monomials in the standard boson algebra.
In particular we will show that
$\Sigma_{\bar{q}}:{\goth H}_{\bar{q}}\to \{W^{AS}_{\bar{q}},\bullet \}$
is an algebra homomorphism which will lead to a phase-space
representation of the deformed Weyl algebra (\ref{wqbar:boson}) analogous
to that of the Weyl algebra (\ref{boson}) provided by (\ref{aso:prod}).
By explicit computation using (\ref{q_ladder}) and (\ref{bigsigma:barq}) we get

\begin{eqnarray}
\sigma_{\bar{q}}({A^+}^n\qs A^m)&=& \sum_{i=0} C(n,m,i)\
\Omega^{\bar{q}}_{n+i\ m+i}            \nonumber \\
C(n,m,i)&=&
\left[{ F_{\bar{q}}!(n+i)\ F_{\bar{q}}!(m+i)\over F_{\bar{q}}!(i)\
F_{\bar{q}}!(i)} \right]^{1/2}. \label{smallsigma:barq}
\end{eqnarray}

\noindent
This transformation can be inverted using

\begin{eqnarray*}
\sigma^{-1}_{\bar{q}}(\Omega^{\bar{q}}_{n m})&=& \sum_{i=0} D(n,m,i)\
{A^+}^{n+i} \qs A^{m+i}
\end{eqnarray*}

\noindent
where the coefficients $D(n,m,i)$ are computed recursively according to

\begin{eqnarray*}
D(n,m,0)&=&C(n,m,0)^{-1}      \\
D(n,m,k)&=&-C(n+k,m+k,0)^{-1}\ \sum_{i=0}^{k-1} C(n+i,n+i,k-i) D(n,m,i).
\end{eqnarray*}

\noindent
$F_{\bar{q}}!(n)\ne 0$ for $0\le n \in {\Bbb Z}$
so $C(n,m,k)$ never vanishes and the coefficients $D(n,m,k)$ are
well-defined. The product
$\bullet : W^{AS}_{\bar{q}}\times W^{AS}_{\bar{q}} \to W^{AS}_{\bar{q}}$
intertwines the deformed Moyal product on $\Omega_{\bar{q}}$ and can easily
be computed using (\ref{algebra_of_states}) as follows

\begin{eqnarray}
w_1 \bullet w_2 &=& \sigma^{-1}_{\bar{q}} (\sigma_{\bar{q}}(w_1) \qs
\sigma_{\bar{q}}(w_2))            \label{qaso:prod}
\end{eqnarray}

\noindent
for $w_1, w_2 \in W^{AS}_{\bar{q}}$. It is a simple matter to show that
$\sigma^{-1}_{\bar{q}}\circ  \sigma_{\bar{q}}$ is the identity on
$W_{\bar{q}}$, and that $\sigma_{\bar{q}}\circ  \sigma^{-1}_{\bar{q}}$
is the identity on $\Omega_{\bar{q}}$.
This establishes the equivalence of $\{W^{AS}_{\bar{q}},\bullet\}$ and
$\{\Omega_{\bar{q}},*\}$ as phase-space representations of $W_{\bar{q}}$
or ${\goth H}_{\bar{q}}$.

Suzuki has considered $*$-products associated with the q-boson which deform
the Moyal $*$-product \cite{suzuki}. By contrast
we have provided $\bullet$-products (\ref{qaso:prod})
which deform the $\bullet$-product of $W^{AS}$ given by (\ref{aso:prod}).
Our construction does not require the introduction on phase-space
of internal degrees of freedom. By constructing
the left-right eigenstate basis for the deformed ring we saw explictly
in \S~\ref{extension} that the Moyal product is not really deformed at all.

\subsection{Infinite matrix representations of ${\goth H}_{\bar{q}}$ and
$W_{\bar{q}}$}\label{infinite_matrices}

Now we define a family of maps
$\pi_{\bar{q}}:\Omega_{\bar{q}} \to M_\infty$, where $M_\infty$ is
the ring of infinite matrices
over the complex numbers. These maps are $\qs$-algebra isomorphisms
of $\{\Omega_{\bar{q}},*\}$ based on (\ref{algebra_of_states}).
We will use the unscripted $\pi:\Omega \to M_\infty$ to denote the action
on the standard boson algebra. Each  $\pi_{\bar{q}}$ is
linear and is defined by its action on $\Omega^{\bar{q}}_{n m}$ according to

\begin{eqnarray}
\pi_{\bar{q}} (\Omega^{\bar{q}}_{n m}) &=& e(n,m)\label{pi:barq}
\end{eqnarray}

\noindent
where $e(n,m)$ is the infinite matrix with zero everywhere except for
the number one in row $n$ and column $m$. The action of
$\pi_{\bar{q}}$ or $\pi_{\bar{q}}\circ \Sigma_{\bar{q}}$ provides
infinite matrix representations of elements of $W_{\bar{q}}$ or
${\goth H}_{\bar{q}}$ respectively in $M_\infty$.

\subsection{Bosonisation of ${\goth H}_{\bar{q}}$ and
$W_{\bar{q}}$}\label{quantum_models}
It is a simple matter to construct bosonic representations of
${\goth H}_{\bar{q}}$ or $W_{\bar{q}}$ either in the space of left-right
eigenstates
$\{\Omega,*\}$ or in $\{W^{AS},\bullet\}$.
In the first case this can be achieved by
the action of $\pi^{-1} \circ \pi_{\bar{q}} \circ \Sigma_{\bar{q}}$
on $u \in {\goth H}_{\bar{q}}$. Applying this to the generators of
${\goth H}_{\bar{q}}$ yields

\begin{eqnarray*}
E&\to&\sum_{n=0}^\infty \Omega_{n n}\\
A^+&\to& \sum_{n=0}^\infty F_{\bar{q}}(n+1)^{1/2}\ \Omega_{n+1\ n}\\
A&\to&\sum_{n=0}^\infty F_{\bar{q}}(n+1)^{1/2}\ \Omega_{n\ n+1}\\
N&\to&\sum_{n=0}^\infty n\ \Omega_{n n}\\
H&\to&{1\over 2} \sum_{n=0}^\infty (F_{\bar{q}}(n+1)+F_{\bar{q}}(n))\
  \Omega_{n n}.
\end{eqnarray*}

A bosonisation of the $\bar{q}$-oscillator can immediately be provided
by

\begin{eqnarray}
A&=& K_{\bar{q}}(N+1) *a \nonumber \\
A^+&=& a^+ *K_{\bar{q}}(N+1) \nonumber \\
H_{\bar{q}}&=&{1\over 2}\left( F_{\bar{q}}(N+1)+F_{\bar{q}}(N) \right)
                                        \label{qbar:model}
\end{eqnarray}

\noindent
where $K_{\bar{q}}(N)= F_{\bar{q}}(N+1)^{1/2} (\hbar(N+1))^{-1/2}$
with $\hbar N = a^+ * a$. There are no ordering problems in
these expressions, but we have to check that
objects such as $K_{\bar{q}}(N)$ are well defined.
{}From {\S}~\ref{extension} we know that
$F_{\bar{q}}$ satisfies $0 < F_{\bar{q}}(n+1)< \infty$ for all
$0\le n \in {\Bbb Z}$. The action of $K_{\bar{q}}(N)$ on $\Omega_{nm}$ is
given by $K_{\bar{q}}(n)= F_{\bar{q}}(n+1)^{1/2} (\hbar(n+1))^{-1/2}$.
The transformation (\ref{qbar:model}) is a well defined
invertible map from $\{ \Omega_{\bar{q}}, *\}$ to $\{ \Omega, * \}$.

Now we return to the main issue of understanding the
non-unitary equivalence of the $\bar{q}$-bosons.

\section{Non-linear representations of boson algebras.}\label{nonlin}
A non-linear representation theory of Lie groups and Lie algebras
developed by Flato, Pinczon and Simon \cite{flato_pinczon_simon}
introduces a concept of non-linear equivalence which is appropriate in
this context. Their work refers specifically
to representations of Lie groups or Lie algebras and so they focus on
transformations which are Lie structure homomorphisms. We are
not dealing with Lie algebras, but with categories of
deformed Weyl or Heisenberg Lie algebras. We therefore consider equivalence
transformations which preserve these categories. In what follows we stick
closely to the intuitive framework provided by Flato et al.

Consider Banach algebras $U$ and $V$. Let
$L_n(U,V)$ denote the $n$-linear maps
{}from $U\times \cdots \times U \to V$. The space of
non-linear transformations $S:U\to V$
consists of maps $s\in S$ such that
$s(u)=\sum_{n=1}^{\infty} f_n(u)$ where $u\in U$ and each $f_n\in L_n(U,V)$.
Algebras or representations of algebras which are related by
invertible transformations of this type will be considered equivalent up
to non-linear transformation.
In particular we would like to consider
$\{\Omega_{\bar{q}},*\}$ or $\{W_{\bar{q}},\bullet\}$ for
$U$ and $V$, in which case the Banach space structure will be derived
{}from the weak topology associated with the Hilbert space structure
(\ref{ip:qboson}).

Let us take two deformed Weyl algebras $W_A, W_B \in \Bbb W$
with generators  $A,A^+,N$ and $B,B^+,N$ as in
{\S}~\ref{boson_algebras}. {}From (\ref{qbar:model}) we can equate

\begin{eqnarray*}
A*A^+=F_{A}(N+1)&& B*B^+=F_{B}(N+1).
\end{eqnarray*}

\noindent
Now we consider the following non-linear transformation

\begin{eqnarray}
A&=&K^A_B(N)*B      \label{nonlin:trans} \\
A^+&=&B^+ *K^A_B(N)  \nonumber
\end{eqnarray}

\noindent
where

\begin{eqnarray*}
K^A_B(N)&=&\left[ F_{\bar{q}^A}(N+1)\over F_{\bar{q}^B}(N+1)\right]^{1\over 2}.
\end{eqnarray*}

\noindent
$K^A_B(N)$ satisfies $0<K^A_B(n)<\infty$
for all $0\le n\in {\Bbb Z}$ so this
transformation leaves $N$ invariant, and (\ref{nonlin:trans})
provides an invertible map {}from $W_A$ to $W_B$.
Phase-space realisations of $W_A, W_B\in \Bbb W$ are therefore equivalent up to
transformations of this form. These transformations are not unitary
since they do not preserve the commutation relations (\ref{wqbar:boson}).
We expect such transformations to exist because of the rigidity of
the Weyl algebra \cite{weyl:rigidity}.
This cohomological result guarantees the existence
of a formal transformation {}from one structure to another, but it must be
verified in each case that this formal transformation is well defined
and invertible. An example where this is so, is the standard Fock-space
representation of the $q$-deformed boson for $q\in {\Bbb R}^+$.
However if $q$ lies at a root of unity the transformation
(\ref{nonlin:trans}) cannot be inverted since $K^A_B(n)$ vanishes for
some positive integer $n$. In this case the equivalence is only a formal
one and two algebras will differ in a fundamental way.

Although we do not elaborate in this paper,
it is straightforward to extend this notion of equivalence
to deformations of the Heisenberg Lie algebra ${\goth H}_{\bar{q}}$
(\ref{hqbar:boson}).

\section{Coproduct structures for ${\goth H}_{\bar{q}}$ and $W_{\bar{q}}$}
\label{coproduct}
Any Lie algebra can be provided with a Hopf algebra structure
\cite{takhtajan}. This amounts to providing a pair of
linear algebra homomorphisms called the co-product
$\Delta: A\to A\otimes A$, and the
evaluation map $\epsilon : A\to C$, as well as a linear algebra
anti-homomorphism called the  antipode  $S:A\to A$. $\Delta$ and
$\epsilon$ are required to satisfy a number of consistency conditions
\cite{takhtajan}. It is a simple matter to provide
$U({\goth h}(1))$ with a Hopf-algebra structure \cite{sun}.
This is achieved by defining the action of the co-product $\Delta$,
the evaluation map $\epsilon$ and
the antipode $S$ on the generators of the Lie algebra
and on the unit as follows

\begin{eqnarray*}
\begin{array}{rclrclrcl}
\Delta (1)&=& 1\otimes 1   & \epsilon (1) &=& 1 & S(1) &=&-1 \\
\Delta (E)&=&E\otimes 1 + 1\otimes E &\epsilon (E)&=&0 & S(E)&=&-E\\
\Delta (A)&=& A\otimes 1 + 1\otimes A &\epsilon (A)&=&0 & S(A)&=&-A\\
\Delta (A^+) &=& A^+ \otimes 1 + 1\otimes A^+&\epsilon (A^+)&=&0
                                                       & S(A^+)&=&-A^+.\\
\end{array}&
\end{eqnarray*}

\noindent
One could try to include the number operator $N$ in the canonical way
by assigning $\Delta(N)=N\otimes 1 + 1\otimes N.$
However this is not compatible with
the definition $\hbar N=A^+*A$ since
$\hbar \Delta(N)\ne \Delta(a^+)\cdot \Delta(a)$. To ensure compatability
we must make the following non-canonical assignment

\begin{eqnarray*}
\Delta(N)&=&N\otimes 1 + 1\otimes N + A\otimes A^+ + A^+ \otimes A.
\end{eqnarray*}

\noindent
In either case the evaluation map and the antipode will act according to
$\epsilon(N)=0$ and $S(N)=-N$. The definitions can now be extended
from the generators $A$, $A^+$, $N$ and $E$ to the entire ring using the
homomorphism properties $\Delta( x\cdot y)=\Delta(x)\cdot \Delta(y)$,
$\epsilon( x\cdot y)=\epsilon(x)\cdot \epsilon(y)$ and the
anti-homomorphism property $S( x\cdot y)=S(y)\cdot S(x)$.

One can provide the Weyl-ring $W$ with a co-product
which is an homomorphism of the associative product as follows

\begin{eqnarray}
\begin{array}{rclrclrcl}
\Delta (1)&=& 1\otimes 1 \\
\Delta (A)&=& (A\otimes 1 + 1\otimes A)/\sqrt{2}\\
\Delta (A^+) &=&(A^+ \otimes 1 + 1\otimes A^+)/\sqrt{2}.\\
\end{array}                                        \label{hopf_weyl}
\end{eqnarray}

\noindent
This can then be extended to the number operator
using $\hbar\ \Delta(N)=\Delta(A^+)\cdot \Delta(A)$, and then to the
whole of $W$ in the natural way.

The only co-unit which is homomorphic
to the associative product on $W$ is the zero map. This follows since
$[\epsilon(a),\epsilon(a^+)]=0=\hbar \epsilon(1)$, so that
$\epsilon(A)=\epsilon(1\cdot A)=\epsilon(1)\epsilon(A)=0$
and similarly $\epsilon(A^+)=0$. As a consequence the Weyl
algebra does not have a bi-algebra or a Hopf-algebra
structure although it does have a co-algebra structure which
is an associative algebra homomorphism. It is possible to use the non-linear
non-canonical equivalence transformations in (\ref{nonlin:trans})
to provide each of the $\bar{q}$-deformed Weyl algebras with a
co-algebra structure as follows

\begin{eqnarray}
\Delta(1)&=& 1\otimes 1                                 \label{q_coproduct}\\
\Delta(A)&=& \Delta\dot K(N)\left[
\left(K(N)^{-1}A\right)\otimes 1+
1\otimes \left(K(N)^{-1}A\right)\right]/\sqrt{2}                 \nonumber \\
\Delta(A^+)&=& \left[ \left(A \otimes K(N)\right)\otimes 1+
1\otimes \left(A^+ K^{-1}(N)\right) \right]
\Delta\cdot K(N)/\sqrt{2}\\
\Delta( N )&=& N\otimes 1 + 1\otimes N                           \nonumber \\
&+&\left(K(N)^{-1} A\right)\otimes \left(A^+ K(N)^{-1}\right)
+\left(A^+ K(N)^{-1}\right)\otimes \left(K(N)^{-1} A\right).     \nonumber
\end{eqnarray}

\noindent
This smoothly deforms the natural co-product structure
of the Weyl algebra (\ref{hopf_weyl})
in such as way that it respects the identity
$\hbar N = A^+ A$ in the undeformed limit. In exactly the same way
one can deform the Hopf-algebra structure of the Heisenberg
Lie algebra whenever the transformation function $K(N,E)$ corresponding to
(\ref{nonlin:trans}) is known.

\section{The classical limit of deformed quantum dynamics}
\label{classical_limit}
The quantum phase-space dynamics of the  $\bar{q}$--deformed oscillator
is easy to study. It suffices to take
$\omega \in L^2({\Bbb R}^2, dq\wedge dp)$ and construct {}from this the
corresponding
density $\rho=\omega \qs \bar{\omega}=\sum_{nm} c_{nm}\
\Omega_{n m}$. Time evolution of this density will be
given by

\begin{eqnarray*}
\rho_t(q,p)&=&\exp\left(-\imath{t\over\hbar} H\qs \right)\qs
      \rho_0(q,p)\qs \exp\left(+\imath{t\over\hbar} H\qs \right) \\
&=&\sum_{n,m=0}c_{nm}\ \exp(\imath t\ \omega_{n m} )\
               \Omega_{nm}(q,p),
\end{eqnarray*}

\noindent
where $\omega_{n m}=( F_{\bar{q}}(m+1)+F_{\bar{q}}(m)-
F_{\bar{q}}(n+1)-F_{\bar{q}}(n))/2$. The expectation of observables
will then be given by

$$ \langle O \rangle_t = (2 \pi \hbar)^{-1}\
\int_{{\Bbb R}^2} O(q,p)\ \rho_t (q,p)\ dq\wedge dp.$$

\noindent
The non-linear transformation (\ref{nonlin:trans})
which relates the algebra of the standard oscillator to that of the
$\bar{q}$-deformed oscillator induces the following transformation
in the defining relations of the deformed Weyl algebra

\begin{eqnarray}
\left.
\begin{array}{rcl}
 {[} N, A ]  &=& -A            \\
 {[}N, A^+ ] &=& A^+           \\
{[} A, A^+ ] &=& f_{\bar{q}}(N)\\
\end{array}
\right\} & \rightarrow &
\left\{
\begin{array}{rcl}
{[}N,a ]&=& -a                 \\
{[}N,a^+ ]&=& a^+              \\
{[}a, a^+ ]&=& \hbar     \\
\end{array}                                              \label{nt_1}
\right.
\end{eqnarray}

\noindent
where

\begin{eqnarray}
f_{\bar{q}}(N)&=&F_{\bar{q}}(N+1)-F_{\bar{q}}(N).         \label{f:def}
\end{eqnarray}

\noindent
At the same time the Hamiltonian is transformed according to

\begin{eqnarray*}
H(A,A^+)={1\over 2}(A^+ A + A A^+) &\rightarrow& H(a,a^+)={1\over 2}
(F_{\bar{q}}(N+1)+F_{\bar{q}}(N))   \label{nt_2}
\end{eqnarray*}

\noindent
where $\hbar N = a^+ a$. To understand the geometrical meaning of these
transformations we need to consider the classical limit of all
of these objects. To start with we put $N= H_0 /\hbar-1/2$ where
$H_0=(a a^+ + a^+ a)/2$ is the Hamiltonian of the
simple harmonic oscillator. Equation (\ref{qbar:model}) becomes

$$ H_{\bar{q}}={1\over 2}\left[
F_{\bar{q}}\left({H_0*\over \hbar}+{1\over 2}\right)+
F_{\bar{q}}\left({H_0*\over \hbar}-{1\over 2}\right)
\right].
$$

\noindent
The limit $\hbar\to 0$ does not make any sense
unless the coefficients of $F_{\bar{q}}$ scale appropriately
with $\hbar$. This is achieved if we replace

\begin{eqnarray}
q_n&\to& q_n\hbar^n           \label{q_subs}
\end{eqnarray}

\noindent
in $F_{\bar{q}}$. The quantum Hamiltonian now becomes

$$ H_{\bar{q}}={1\over 2}\left[
F_{\bar{q}}\left(H_0* +{\hbar \over 2}\right)+
F_{\bar{q}}\left(H_0* -{\hbar \over 2}\right)
\right]
$$

\noindent
and the classical limit of the $\bar{q}$-deformed oscillator is
given by

\begin{eqnarray}
H^{C}_{\bar{q}}=F_{\bar{q}}(H_0)&=& q_0+q_1 H_0 + q_2 H^2_0 + \ldots
                              \label{cl:ham}
\end{eqnarray}

\noindent
{} From (\ref{nt_1}) and (\ref{f:def}) the commutator of $A$ and $A^+$
in (\ref{wqbar:boson}) becomes

\begin{eqnarray*}
A\qs A^+ - A^+ \qs A &=&F_{\bar{q}}\left(H_0* +{\hbar \over 2}\right)-
F_{\bar{q}}\left(H_0* -{\hbar \over 2}\right) \\
                        &=& F^{\prime}_{\bar{q}}(H_0 )\hbar  + O(\hbar^3).
\end{eqnarray*}

\noindent
In the classical limit this corresponds to the following
non-canonical Poisson bracket

\begin{eqnarray}
\{ f, g \}_{\bar{q}} &=&  F^{\prime}_{\bar{q}}(H_0)\{f, g \}. \label{new_pb}
\end{eqnarray}

\section{Relationship to method of normal forms}\label{normal_forms}

In classical Hamiltonian dynamics on the plane an Hamiltonian $H(q,p)$
in  the  neighbourhood of  a   point  of  stable   equilibrium  can be
transformed by a series of canonical  transformations to a co-ordinate
system in which it takes the  form $H(p^2+q^2)$. These transformations
preserve  the Poisson bracket  and their quantum analogues are unitary
transformations in the quantum version of  the method of  normal forms
\cite{crehan_nforms}.  Effectively they smooth the invariant curves of
the flow into Euclidean circles. Expansion in normal forms is an
important   technique to   study  the non-linear    dynamics of a
system in the neighbourhood of a critical point
\cite[Chapter 4, eqn 4.96]{normal_form}.  The truncated  expansion  can provide
a useful
integrable approximation to a non-integrable system.  This
provides the starting  point  for an important semi-classical approach
to  the  quantisation of a range of non-integrable  Hamiltonian
systems such as small molecules.

$H^{C}_{\bar{q}}$ (\ref{cl:ham}) is a power series expansion  in terms
of the classical Hamiltonian of the simple  harmonic oscillator. It is
the   expansion  in  normal   forms  of  a  one-dimensional  classical
Hamiltonian    system  in the   neighbourhood   of a  point of  stable
equilibrium. Up to a simple  linear transformation the coefficients of
powers  of $H_0$ in  $H^{C}_{\bar{q}}$  in (\ref{cl:ham}) are just the
$\bar{q}$-deformation   parameters   which appeared   originally    in
(\ref{wqbar:boson}). Hamiltonian  systems of this form  are  very like
linear oscillators  in   the sense that the   orbits are  confined  to
circles about the origin. Non-linearity  is exhibited  in the fact
that the fundamental frequency of each orbit  depends on the circle to
which it is confined.

It is possible to  make a further change of  co-ordinates so that  a
Hamiltonian  of the  form $H(p^2+q^2)$  is transformed into  the  even
simpler form $(p^2+q^2)/2$.  This last transformation is not canonical
since the Poisson  bracket changes  form. It is  in fact the classical
counterpart of  (\ref{nonlin:trans})  which relates  the deformed Weyl
algebras   (\ref{algebra_of_states})  to  the  Weyl   algebra   itself
(\ref{boson}).  The   inverse of (\ref{nonlin:trans})   which  in  the
classical limit would transform $H(q,p)$ into the  form $(p^2+q^2)/2$,
seems   closely   related  to  the  factorisation    method for linear
differential  equations  on   the   line due  to    Infeld    and Hull
\cite{infeld_and_hull}.

\section{Quantisation of non-canonical Poisson structures}\label{quantisation}
In  previous  work  on  the existence   of Hamiltonian  structures for
arbitrary vector fields
\cite{crehan_2}, \cite{crehan_3}, one of us raised the question as to whether
the equivalence between the infinite number of Hamiltonian structures
of a given classical Hamiltonian system extends to its quantisation.
The simple harmonic oscillator is not just Hamiltonian with
Hamiltonian function $H_0 =a a^+$ with respect to the Poisson bracket
$\{a,a^+\}=1$, but also with Hamiltonian function $H=F(H_0)$ with
respect to the Poisson bracket  $\{a,a^+\}=(F^\prime (H_0))^{-1}$.
It has an infinite number of non-canonical Hamiltonian
structures. A natural question arises
as to whether there exists a theory of
quantisation which consistently applies to these
classically equivalent formulations. On the basis of
(\ref{f:def}), (\ref{q_subs}) and (\ref{new_pb}) we can recast our
understanding of the classical limit of the deformed commutation
relations

\begin{eqnarray*}
\{a,a^+\}=F^{\prime}_{\bar{q}}(H_0)&\leftarrow&
{[}a,a^+]= F_{\bar{q}}(H_0+\hbar/2)-F_{\bar{q}}(H_0-\hbar/2)
\end{eqnarray*}

\noindent
as a quantisation rule which applies to these non-canonical Poisson
structures

\begin{eqnarray}
\{a,a^+\}=\Theta_{\bar{q}}(H_0)\rightarrow {[}a,a^+]&=&
\int^{H_0 +\hbar/2 }_{H_0-\hbar/2}\Theta_{\bar{q}}(x)\ dx   \label{ncqs}\\
&=& \hbar\ \Theta_{\bar{q}}(H_0) + {\hbar^3\over 24}\
\Theta_{\bar{q}}^{\prime \prime}(H_0) +\ldots
                             \nonumber
\end{eqnarray}

\noindent
Looking at the deformed Weyl algebras (\ref{wqbar:boson})
{}from this point of view
we see that the fundamental commutation relation is not only
a deformation of the classical Poisson bracket, but
also a coarse graining of the classical Poisson bracket whose
coarseness is controlled by Planck's constant. The quantum
commutator is also the spectral--gap operator for the energy
levels of a non-linear osillator.

\section{Discussion}
{}From cohomological arguments we know that formally there are
no non-trivial associative-deformations of the Weyl-algebra
\cite{weyl:rigidity}. This raises the problem of how to understand
quantum deformations of the Weyl algebra such as the $q$-- and
$qp$--deformed boson algebras.

To tackle this we introduce a family of
deformed Weyl and Heisenberg Lie algebras which
depend on an infinite number of independent deformation parameters.
Using phase-space quantum mechanical methods
we studied phase-space realisations of these algebras using ASO monomials
and bimodular bases of deformed-oscillator eigenstates.
In the monomial bases the intertwiner of the deformed
Moyal product was given by a $\bullet$-product (\ref{qaso:prod}) which
deforms that associated with the standard boson (\ref{aso:prod}).
We provided an algorithm for transforming to the eigenstate bases.
In the eigenstate bases the deformed Moyal--products are immediately seen to
be isomorphic to that of the standard oscillator,  and
to the algebra of infinite matrices. These results are
best sumarised in the following diagram and by equations
(\ref{algebra_of_states}) and (\ref{ip:qboson}). The diagram relates
the quantum deformed Heisenberg Lie algebra ${\goth H}_{\bar{q}}$, to the
phase-space representation of the deformed Weyl algebra $W_{\bar{q}}$
given by $\{W^{AS}_{\bar{q}},\bullet\}$ and $\{\Omega^{\bar{q}},*\}$,
to $M_{\infty}$ the algebra of infinite matrices, and to representations
of the Weyl algebra given by $\{W^{AS},\bullet\}$ and $\{\Omega,*\}$.
The vertical arrows are algebra homomorphisms, whereas the
horizontal arrows are isomorphisms of phase-space representations
of the deformed Weyl algebras (\ref{wqbar:boson}) which satisfy the
vacuum postulate (\ref{qbar:vacuum}).

\begin{eqnarray}
&
\begin{array}{rcccccccc}
&&{\goth H}_{\bar{q}}&&&&&& \\
&&&&&&&&\\
&&\downarrow&&&&&& \\
&&&&&&&&\\
&&W_{\bar{q}}&&&&&&\\
&&&&&&&&\\
&&\downarrow&&&&&& \\
&&&&&&&&\\
\{W^{AS}_{\bar{q}},\bullet\}&
{\rightleftharpoons}&   \{\Omega_{\bar{q}},*\}&
{\rightleftharpoons}&           M_\infty&
{\rightleftharpoons}&            \{\Omega,*\} &
{\rightleftharpoons}&\{ W^{AS},\bullet\}\\
\end{array}
&\label{main_theorem}
\end{eqnarray}

\noindent
This enables us to construct phase-space quantum mechanical
models of the deformed Weyl algebras and their associated oscillators.
The equivalence  of all realisations of (\ref{wqbar:boson})
satisfying the vacuum postulate (\ref{qbar:vacuum}), is not a
unitary equivalence since it does not preserve the form of the
commutation relations. These algebras are equivalent
in the sense that they are related by non-linear non-unitary
transformations (\ref{nonlin:trans}). This allows
us to provide the deformed Weyl-algebra with a
co-product structure (\ref{q_coproduct}) which is a
deformed Weyl algebra homomorphism and which smoothly deforms the
co-product structure of the Weyl algebra (\ref{hopf_weyl}).

Finally we considered the classical limit of
phase-space dynamics of the deformed
oscillators. We found that
the deformation parameters have a natural interpretation
as coefficients of non-linear terms in a normal forms expansion
of the deformed oscillator (\ref{cl:ham}).
The situation is best summed up by
considering a fully integrable classical system whose
Hamiltonian is an entire complex function of the harmonic oscillator
$F_{\bar{q}}(H_0)=q_0 + q_1 H_0 + q_2 H_0^2 +\ldots$ with respect to
the Poisson bracket $\{ A, A^+ \} = 1$. The parameters
$q_n$ for $2\le n \in {\Bbb Z}$ parameterise non-linear deformations of
the harmonic oscillator. The quantum Hamiltonian is given by
$(F_{\bar{q}}(H_0+\hbar/2)+ F_{\bar{q}}(H_0-\hbar/2))/2$.
It can be transformed into
$H_0=(A A^+ + A^+ A)/2$ by a non-linear non-canonical transformation of
co-ordinates. Under this transformation the fundamental
quantum commutator $[A,A^+]=\hbar$ becomes
$[A,A^+]=(F_{\bar{q}}(H_0+\hbar/2)-F_{\bar{q}}(H_0-\hbar/2)$ which we recognise
as a deformation of the Weyl algebra of the form (\ref{wqbar:boson}).
In the classical limit the transformed Hamiltonian is $A A^+$ and the Poisson
bracket is given by $\{ A, A^+ \} = F^{\prime}_{\bar{q}}(H_0)$.
Eventually we are lead to consider the deformed Weyl algebras
(\ref{wqbar:boson}) as providing fundamental
commutation relations (\ref{ncqs}) which quantise the non-canonical
classical Poisson brackets given in (\ref{new_pb}). Other
families of deformations of the Weyl-algebra can easily be
constructed. This would lead to a quantisation of the infinite family of
non-canonical Hamiltonian structures associated with more general
one dimensional systems.

Classically we can consider the family of quantum--deformations of the
Weyl algebra considered in \S~\ref{boson_algebras},
as being equivalent to non-linear deformations of the classical
oscillator (\ref{cl:ham}), or non-canonical
transformations of the related Poisson structure
(\ref{new_pb}). Classically the deformed Poisson algebra is equivalent
to the undeformed algebra as long as $F_{\bar{q}}(z)$ is an invertible map
on ${\Bbb R}^+$. By contrast phase-space realisations of the
related quantum algebras are equivalent as
long as $0 < F_{\bar{q}}(n+1) < \infty$. The $q$-deformed
oscillator for $q$ a root of unity is of considerable interest in physics
\cite{root_of_unity}, and provides a good example of how formal or local
equivalence guaranteed by cohomology does not
extend to global equivalence between the deformed and undeformed systems.

Realisations  of  the $q$-boson exist  in  terms of   standard boson
eigenstates on configuration  space \cite{kulish_2}, \cite{kibler}, on
complexified  configuration   space \cite{kempf},   or  in  terms   of
q-difference operators acting  on the functions of configuration space
\cite{atakishiyev}.  We think that phase-space quantum mechanical
techniques can lead to a
clarification of the similarities and differences between the standard
boson  and  the  $q$-boson  which  arise  in  the  configuration-space
formalism of standard quantum mechanics. We have considered deformations
of the Weyl algebra which satisfy the vacuum postulate. This
is the smallest category within which to deform the Weyl algebra.
Rideau has shown that there
exist infinite  matrix representations  of the $q$-boson for which the
vacuum  postulate  is  not  satisfied \cite{rideau}.  In view of the
existence of these other representations it is of
interest to extend our results to the larger category of
representations which do not necessarily satisfy the  vacuum postulate.

It   should be straightforward   to  provide  quantum  deformations of
semi-simple  Lie   algebras such as   ${\goth {su}}(2)$
using  the Jordan-Schwinger
construction.  It  will be of interest  to  extend the results presented
in this paper to the higher dimensional Weyl algebras $W(n)$.
This should provide some new perspectives on the role  of  the
$R$-matrix  in  quantum  group constructions \cite{r_matrix}.
Quantum groups  now find
application in   nuclear  physics, molecular dynamics, and  mesoscopic
systems as well as in  the quantum field  theory of integrable  models
\cite{kibler}, \cite{kulish_1}, \cite{takhtajan}. We expect  that our results
on   deformed  boson  algebras  provide
a physically intuitive way to understand quantum deformation, which
will be of use  in
generalising applications based on  standard $q$-- and $qp$--deformed
boson algebras.

\section{Acknowledgement}
P.Crehan would like to thank Michio Jimbo for his hospitality at
Kyoto University. T.G.Ho would like to thank Prof. Araki of the
Research Institute for
Mathematical Sciences. We would both like to thank Moshe Flato and
Daniel Sternheimer for their many valuable comments and criticisms as well as
for bringing some important references to our attention.

\end{document}